\documentclass[11pt,a4paper]{article}

\usepackage{amssymb}
\usepackage{amsmath}
\usepackage{amsthm}
\usepackage{amscd}
\usepackage[all]{xy}

\newtheorem{Thm}{Theorem}[section]

\newtheorem{Lem}[Thm]{Lemma}
\newtheorem{Prop}[Thm]{Proposition}

\theoremstyle{remark}

\newtheorem*{pf}{\rm\textbf{Proof}}

\newtheorem*{ac}{\rm\textbf{Acknowledgement}}

\begin{document}

\title{\Large \bf  Non-simply-laced Clusters of Finite Type via Frobenius Morphism}

\author{\large Dong YANG\\
{\small Department of Mathematical Sciences,}\\
{\small Tsinghua University,}\\
{\small Beijing100084, P.R.China.}\\
{\small yangdong98@mails.tsinghua.edu.cn.} }

\date{}

\maketitle

\begin{quote} \abstract{ By showing the compatibility of folding almost positive roots
 and folding cluster categories, we prove that there is a one-to-one correspondence between seeds and tilting seeds in non-simply-laced finite cases.}
\end{quote}
\vskip30pt

\section{Introduction}

Introduced by Fomin and Zelevinsky~\cite{FZ}, cluster algebras are
a family of commutative algebras generated by the so-called
cluster variables. In~\cite{FZ2} the authors classify all cluster
algebras of finite type, and construct for all finite types a
one-to-one correspondence between almost positive roots and
cluster variables. It is pointed out that this classification
coincides with the classification of Cartan matrices.
Then~\cite{MRZ}~\cite{BMRRT} (and independently~\cite{CCS})
developed the representation theoretical approach to study cluster
algebras.~\cite{BMR}~\cite{CC}~\cite{CK}~\cite{CK2}~\cite{Hub2}~\cite{Z}
are papers following this approach. In particular, they establish
a one-to-one correspondence between clusters and tilting objects
in the so-called cluster category.

In~\cite{MRZ} the authors claim that non-simply-laced clusters can
be studied by folding the corresponding simply-laced ones. Using
this method Dupont~\cite{D} realizes non-simply-laced cluster
algebras as quotients of simply-laced cluster algebras. Deng and
Du~\cite{DD} establish a link between representations of quivers
over $\overline{\mathbb{F}}_{q}$ (simply-laced case) and
representations of $\mathbb{F}_{q}$-species (non-simply-laced
case). This link is generalized in~\cite{DD2} to the homotopy
category, the bounded derived category and the root category. In
this paper we extend it further to the cluster category. Namely,
we study folding cluster categories. Applying this and the study
of folding almost positive roots to known results on simply-laced
clusters of finite type, we provide a new approach to understand
non-simply-laced clusters of finite type. The main result is

\begin{Thm} For non-simply-laced finite type, there is a one-to-one correspondence between clusters
and tilting objects in the cluster category. Under this
correspondence, mutation of seeds is exactly mutation of tilting
seeds.
\end{Thm}

\begin{ac} The author would like to thank Steffen Koenig and
Department of Mathematics at University of Koeln for their
hospitality during his visit, when this work started. He also
thanks Jie Xiao and Bin Zhu for helpful conversations. The author
acknowledges support by National Natural Science Foundation of
China No.10131010.
\end{ac}

\section{Preliminaries}
We first introduce some notations.

For a category $\mathcal{E}$, denote by $ind\mathcal{E}$ the set
of isomorphism classes of indecomposable objects of $\mathcal{E}$.
Depending on the context, it also means a representative set. For
an object $X$ of $\mathcal{E}$, denote by $[X]$ the isomorphism
class of $X$.

Let $A$ be an algebra. By an $A$-module, we mean a finite
dimensional left $A$-module. Denote by $modA$ the category of
$A$-modules.

Let $q$ be a prime power, $\mathbb{F}_{q}$ the finite field with
$q$ elements, and $\overline{\mathbb{F}}_{q}$ a fixed algebraic
closure of $\mathbb{F}_{q}$.

\subsection{Frobenius morphisms}\label{S:frobmor}

In~\cite{DD}~\cite{DD2}~\cite{DD3}, Deng and Du study relation
between the representation theory of an
$\overline{\mathbb{F}}_{q}$-algebra $A$ with a Frobenius
automorphism $F$ and that of the algebra $A^{F}$ of fixed points,
which is an $\mathbb{F}_{q}$-algebra. This subsection is devoted
to a brief introduction to their results.

A {\em Frobenius map} $F$ on a $\overline{\mathbb{F}}_{q}$-space
$V$ is a map from $V$ to itself satisfying

(i) $F(\lambda v)=\lambda^{q}v$, $\forall
\lambda\in\overline{\mathbb{F}}_{q},\ v\in V$.

(ii) $\forall v\in V$, $\exists n\in\mathbb{N}$, s.t.
$F^{n}(v)=v$.

\begin{Lem} {\rm (i)} Let $V$ be a finite dimensional $\overline{\mathbb{F}}_{q}$-space with a
Frobenius map $F_{V}$, then $V^{F}=\{v\in V\ |\ F_{V}(v)=v\}$ is a
finite dimensional $\mathbb{F}_{q}$-space such that
$V=\overline{\mathbb{F}}_{q}\otimes_{\mathbb{F}_{q}}V^{F}$.
Moreover, the $\mathbb{F}_{q}$-dimension of $V^{F}$ equals the
$\overline{\mathbb{F}}_{q}$-dimension of $V$.

{\rm (ii)} Let $V$, $W$ be two finite dimensional
$\overline{\mathbb{F}}_{q}$-spaces with Frobenius maps $F_{V}$ and
$F_{W}$ respectively, then $F:f\mapsto F_{W}\circ f\circ
F_{V}^{-1}$ is a Frobenius map on
$Hom_{\overline{\mathbb{F}}_{q}}(V,W)$, and
$Hom_{\mathbb{F}_{q}}(V^{F},W^{F})=(Hom_{\overline{\mathbb{F}}_{q}}(V,W))^{F}$.

\end{Lem}

In the rest of this subsection we assume $A$ is a finite
dimensional algebra over $\overline{\mathbb{F}}_{q}$ with a fixed
Frobenius morphism. By definition a {\em Frobenius morphism} $F$
on $A$ is a Frobenius map on the vector space $A$ and also an
automorphism of $\mathbb{F}_{q}$-algebras. Note that the algebra
$A^{F}$ of fixed points is a finite dimensional algebra over
$\mathbb{F}_{q}$, and $\overline{\mathbb{F}}_{q}\otimes A^{F}=A$.
Moreover, $A$ is hereditary if and only if $A^{F}$ is hereditary.

For an $A$-module $M$ and $r\in\mathbb{Z}$, define an $A$-module
$M^{[r]}=A\otimes_{F^{r}}M$, where $a\otimes m=1\otimes
F^{-r}(a)m$. If we write $m^{(r)}$ for $1\otimes m$ and write the
new action by $\cdot_{r}$, then explicitly
$a\cdot_{r}m^{(r)}=(F^{-r}(a)m)^{(r)}$. For an $A$-module
homomorphism $f:M\rightarrow N$, we define an $A$-module
homomorphism $f^{[r]}:M^{[r]}\rightarrow N^{[r]}, m^{(r)}\mapsto
(f(m))^{(r)}$. In this way, $()^{[r]}$ is an autoequivalence of
$modA$ and $()^{[r]}()^{[s]}=()^{[r+s]}$. In particular,
$()^{[r]}$ commutes with the Auslander-Reiten translation $\tau$.

The selfequivalence $()^{[r]}$ of the Abelian category $mod A$
extends to $\mathcal{C}(A)$, the category of complexes of
$A$-modules, which commutes with the shift functor $[1]$.
Precisely, for $M^{\cdot}=\{M_{i},d_{i}\}$, we define
$\{M_{i},d_{i}\}^{[r]}=\{M_{i}^{[r]},d_{i}^{[r]}\}$, or
equivalently, ${M^{\cdot}}^{[r]}=A\otimes_{F^{r}}M^{\cdot}$. For
$f^{\cdot}=\{f_{i}\}:M^{\cdot}\rightarrow N^{\cdot}$ a
homomorphism of complexes,
${f^{\cdot}}^{[r]}=\{f_{i}^{[r]}\}:{M^{\cdot}}^{[r]}\rightarrow
{N^{\cdot}}^{[r]}$ is a homomorphism of complexes too. Moreover,
$f$ is homotopic to zero (resp. a quasi-isomorphism) if and only
if so is $f^{[r]}$. Therefore $()^{[r]}$ induces an equivalence of
triangulated categories $()^{[r]}:C(A)\rightarrow C(A)$ where
$C(A)$ is the bounded homotopy category $\mathcal{K}^{b}(A)$ or
the bounded derived category $\mathcal{D}^{b}(A)$.

In the following let $C(A)=modA, \mathcal{C}^{b}(A),
\mathcal{K}^{b}(A),$ or $\mathcal{D}^{b}(A)$. We have
$dim_{\overline{\mathbb{F}}_{q}}Hom_{C(A)}(M,N)=dim_{\overline{\mathbb{F}}_{q}}Hom_{C(A)}(M^{[1]},N^{[1]})$.
A pair $(M,\phi_{M})$ is called an {\em $F$-stable object} in
$C(A)$ if $\phi_{M}:{M}^{[1]}\rightarrow M$ is an isomorphism in
$C(A)$. We also say that $M$ is an $F$-stable object. Denote by
$C^{F}(A)$ the category of $F$-stable objects whose morphisms
$f:(M,\phi_{M})\rightarrow (N,\phi_{N})$ are morphisms
$f:M\rightarrow N$ in $C(A)$ satisfying
$f\circ\phi_{M}=\phi_{N}\circ {f}^{[1]}$. On $C^{F}(A)$ we have
the Auslander-Reiten translation $\tau$ defined by
$\tau(M,\phi_{M})=(\tau_{M},\tau\phi_{M})$.

Let $(M,\phi_{M})$ be an $F$-stable $A$-module, i.e. an $F$-stable
object in $modA$. Then $F_{M}:M\rightarrow M$, $m\mapsto
\phi_{M}(m^{(1)})$ is a Frobenius map on $M$ which is compatible
with the $A$-module structure of $M$, i.e.
$F_{M}(am)=F(a)F_{M}(m)$ for any $a\in A,m\in M$. We see that
$M^{F}$ is an $A^{F}$-module. An morphism
$f:(M,\phi_{M})\rightarrow (N,\phi_{N})$ between $F$-stable
$A$-modules induces a morphism $f^{F}=f|_{M^{F}}:M^{F}\rightarrow
N^{F}$ between $A^{F}$-modules. Thus we have a functor $\Phi$ from
$C^{F}(A)$ to $C(A^{F})$. For an object $(M,\phi_{M})$ in
$C^{F}(A)$ we will denote by $M^{F}$ its image under $\Phi$. We
have

\begin{Thm}\label{T:dd} {\rm (i)} $\Phi$ is an equivalence (of abelian categories and triangulated
categories respectively) from $C^{F}(A)$  to $C(A^{F})$ with
inverse functor
$\overline{\mathbb{F}}_{q}\otimes_{\mathbb{F}_{q}}-$. Moreover,
this equivalence commutes with the Auslander-Reiten translation.

{\rm (ii)} Let $(M,\phi_{M})$ and $(N,\phi_{N})$ be two objects in
$C^{F}(A)$. Then $M^{F}$ and $N^{F}$ are isomorphic in $C(A^{F})$
if and only if $M$ and $N$ are isomorphic in $C(A)$. Moreover,
$Hom_{C(A^{F})}(M^{F},N^{F})=(Hom_{C(A)}(M,N))^{F}$.

\end{Thm}

It is shown that for each object $M$ in $C(A)$ there exists
$r\in\mathbb{N}$ such that $M^{[r]}\cong M$. For such $r$ the
module $M\oplus M^{[1]}\cdots\oplus M^{[r-1]}$ is $F$-stable. If
$r$ is minimal such that $M^{[r]}\cong M$ (such $r$ is called the
$F$-period of $M$) then $\tilde{M}=M\oplus M^{[1]}\cdots\oplus
M^{[r-1]}$ is indecomposable in $C^{F}(A)$.

\subsection{Clusters}

In this subsection we follow~\cite{FZ2}~\cite{FZ3}. Let $\Gamma$
be a valued graph of finite type with vertex set $I$, denote by
$\Phi(\Gamma)$ the root system and by $\Phi_{\geq -1}(\Gamma)$ the
set of almost positive roots (i.e. the set of positive roots and
negative simple roots). A {\em cluster variable} is an element in
$\Phi_{\geq -1}(\Gamma)$.

Let $s_{i}$ be the simple reflection of the Weyl group of
$\Phi(\Gamma)$ corresponding to $i\in I$, and let $\sigma_{i}$ be
the permutation of $\Phi_{\geq -1}(\Gamma)$ defined as follows
\[\sigma_{i}(\alpha)=\begin{cases} \alpha & {\rm\text if\ }\alpha=-\alpha_{j},j\not=i\\
                               s_{i}(\alpha) & {\rm \text otherwise}
                 \end{cases}\]

Let $I=I^{+}\cup I^{-}$ be a partition of the vertex set $I$ into
completely disconnected subsets. Let $\tau_{+}=\prod_{i\in
I^{+}}\sigma_{i}$ and $\tau_{-}=\prod_{i\in I^{-}}\sigma_{i}$.
Define the compatibility degree
$(\alpha||\beta)_{\Gamma}=(\alpha||\beta)$ by
$(-\alpha_{i}||\beta)={\rm \text the\ coefficient\ of\ }\alpha_{i}
{\rm \text\ in\ }\beta$, where $\alpha_{i}$ is a simple root and
$\beta$ is a positive root, and extend $\tau_{\pm}$-invariantly to
$(\ ,\ ):\Phi_{\geq -1}(\Gamma)\times \Phi_{\geq
-1}(\Gamma)\rightarrow \mathbb{N}_{0}$. If $\Gamma$ is
simply-laced then $(\alpha||\beta)=(\beta||\alpha)$.

We say that $\alpha,\beta\in\Phi_{\geq -1}(\Gamma)$ are {\em
compatible} if $(\alpha||\beta)=0$. A {\em cluster} is defined to
be a maximal compatible subset of $\Phi_{\geq -1}(\Gamma)$. The
cardinality of a cluster equals the cardinality of $I$. Two
elements $\alpha,\beta$ ($\alpha\neq \beta$)\ in $\Phi_{\geq
-1}(\Gamma)$ {\em form an exchange pair} if there is a cluster
$\underline{x}$ such that $\alpha\in \underline{x}$ and
$(\underline{x}\backslash \{\alpha\})\cup \{\beta\}$ is again a
cluster.

Let $\underline{x}=\{x_{i}\}_{i\in I}$ a cluster, and $B=(b_{ij})$
be any matrix realization of the graph $\Gamma$, which is an
anti-symmetrizable matrix. We call the pair $(\underline{x},B)$ an
{\em initial seed}. For each $k\in I$, one can define the {\em
mutation} along direction $k$ which is a pair
$(\underline{x}',B')$ where
$\underline{x}'=(\underline{x}\backslash x_{k})\cup \{x_{k}'\}$
with $x_{k}'$ defined via $\underline{x}$ and $B$ (in fact $x_{k}$
and $x_{k}'$ form an exchange pair), and $B'=(b_{ij}')$ is defined
by
\[b_{ij}'=\begin{cases} -b_{ij}  &  {\rm \text\ if\ } i=k {\rm \text\ or\ } j=k\\
                        b_{ij}+\frac{b_{ik}|b_{ij}|+|b_{ik}|b_{kj}}{2}
                        &{\rm\text \ otherwise}
                        \end{cases}\]
We write $(\underline{x}',B')=\mu_{k}(\underline{x},B)$, and
$B'=\mu_{k}B$. We call a pair $(\underline{x}',B')$ a {\em seed}
if $\underline{x}'$ is a cluster and $B'$ is an anti-symmetrizable
matrix such that the pair is an iterated mutation of the initial
seed.

To each anti-symmetrizable matrix $B=(b_{ij})$ with integer
entries we associate a valued quiver $Q_{B}$. Precisely, if
$b_{ij}>0$, then $Q_{B}$ has $ \UseTips\entrymodifiers={++[][]}
\xymatrix @=10mm { i \ar[r]^{(-b_{ji},b_{ij})} & \ j}$. Note that
$Q_{B}$ has no loops or oriented cycles of length $2$. In fact
this defines a bijection between the set of anti-symmetrizable
matrices with integer entries and the set of valued quivers
without loops or oriented cycles of length $2$.  In the following,
we denote by $B_{Q}$ the anti-symmetrizable matrix corresponding
to a valued quiver $Q$.

\subsection{Cluster category}
 Let $A$ be a finite dimensional hereditary algebra. Let $\tau$ be
 the Auslander-Reiten translation of $\mathcal{D}^{b}(A)$ and
 let $S$ denote the shift functor. Following~\cite{BMRRT} we form the orbit category
 $\mathcal{C}_{A}=\mathcal{D}^{b}(A)/\tau^{-1}S$. This is a
 Krull-Schmidt category.
 Precisely, the objects of $\mathcal{C}_{A}$ are exactly the
 objects in $\mathcal{D}^{b}(A)$, and for $X,Y\in\mathcal{C}_{A}$,
 $Hom_{\mathcal{C}_{A}}(X,Y)=\oplus_{i\in\mathbb{Z}}Hom_{\mathcal{D}^{b}(A)}(X,(\tau^{-1}S)^{i}Y)$. It is shown in~\cite{Ke} that $\mathcal{C}_{A}$ is a triangulated
 category with shift functor $S$, and the canonical functor from $\mathcal{D}^{b}(A)$ to $\mathcal{C}_{A}$ is a triangle functor.

Define
 $Ext_{\mathcal{C}_{A}}^{i}(X,Y)=Hom_{\mathcal{C}_{A}}(X,S^{i}Y)$ for $i\in\mathbb{Z}$.
 Then
 $Ext_{\mathcal{C}_{A}}^{1}(X,Y)=Ext_{\mathcal{C}_{A}}^{1}(Y,X)$.
 An object $T$ in $\mathcal{C}_{A}$ is called {\em exceptional} if
 $Ext^{1}_{\mathcal{C}_{A}}(T,T)=0$.
An exceptional object $T$ is called a {\em (basic cluster)-tilting
object} (resp. {\em almost complete tilting object}) if $T$ has
exactly $n$ (resp. $n-1$) pairwise non-isomorphic indecomposable
direct summands, where $n$ is the number of isoclasses of simple
$A$-modules. An exceptional indecomposable object $M$ is called a
{\em complement} of an almost complete tilting object
$\overline{T}$ if $\overline{T}\oplus M$ is a tilting object. An
almost complete tilting object has precisely two complements.
Moreover, for two exceptional indecomposable objects $M$ and
$M^{*}$, there exists an almost complete tilting object
$\overline{T}$ such that $M$ and $M^{*}$ are exactly the two
complements of $\overline{T}$ if and only if
$dim_{D_{M}}Ext^{1}_{\mathcal{C}_{A}}(M,M^{*})=1=dim_{D_{M}^{*}}Ext^{1}_{\mathcal{C}_{A}}(M^{*},M)$
where $D_{M}$ and $D_{M^{*}}$ are the division algebras
$End(M)/radEnd(M)$ and $End(M^{*})/radEnd(M^{*})$ respectively.

Let $T=T_{1}\oplus\cdots\oplus T_{n}$ be a tilting object. The
algebra $End_{\mathcal{C}_{A}}(T)^{op}$ is called a {\em
cluster-tilted algebra}. By~\cite{BMR} Proposition 3.2  the valued
quiver $Q_{T}$ of $End_{\mathcal{C}_{A}}(T)^{op}$ has no loops or
oriented cycles of length $2$. Write $B_{T}=B_{Q_{T}}$. The pair
$(T,B_{T})$ is called a {\em tilting seed}. For $k=1,\cdots,n$,
there exists a unique indecomposable object $T_{k}^{*}$
non-isomorphic to $T_{k}$ such that $T'=T_{1}\oplus\cdots\oplus
T_{k-1}\oplus T_{k}^{*}\oplus T_{k+1}\oplus\cdots\oplus T_{n}$ is
again a tilting object. Therefore $(T',B_{T'})$ is again a tilting
seed, called a {\em mutation} of $(T,B_{T})$ along direction $k$.
We write $(T',B_{T'})=\delta_{k}(T,B_{T})$, and
$B_{T'}=\delta_{k}B_{T}$. Any tilting seed is an iterated mutation
of the tilting seed $(SA,B_{SA})$.

\subsection{Simply-laced clusters of finite type}

Let $Q$ be a Dynkin quiver with underlying graph $\Gamma$ (hence
$\Gamma$ is a simply-laced Dynkin graph of finite type). Denote by
$I$ the set of vertices of $Q$ (and $\Gamma$).

Let $P_{i}$ be the indecomposable projective module corresponding
to the vertex $i$. Then $ind\mathcal{C}_{A}=ind(modA)\cup
\{[SP_{i}]\}_{i=1,\cdots,n}$. We define the dimension vector of
$SP_{i}$ to be $\underline{dim}(SP_{i})=-\alpha_{i}$. By Gabriel's
theorem (see~\cite{G}) taking dimension vector defines a bijective
correspondence between the set $\{SP_{i}\ |\ i\in I\}\cup ind\
mod(\overline{\mathbb{F}}_{q}Q)$ and $\Phi_{\geq -1}(\Gamma)$,
namely, between $ind\mathcal{C}(Q)$ and the set of cluster
variables.

\begin{Thm}\label{T:BMRRT}(~\cite{BMRRT})
{\rm (i)} For $\alpha,\beta\in\Phi_{\geq -1}(\Gamma)$, let
$M_{\alpha},M_{\beta}$ be the corresponding indecomposable object
in $\mathcal{C}(Q)$, then
$(\alpha||\beta)_{\Gamma}=dim_{\overline{\mathbb{F}}_{q}}Ext^{1}_{\mathcal{C}(Q)}(M_{\alpha},M_{\beta})$.

{\rm (ii)} The bijective correspondence induces a bijective
correspondence between the isoclasses of basic (cluster-) tilting
objects in $\mathcal{C}(Q)$ and clusters.

{\rm (iii)} Two elements $\alpha,\beta$ of $\Phi_{\geq
-1}(\Gamma)$ form an exchange pair if and only if their
compatibility degree is $1$.
\end{Thm}

Theorem~\ref{T:BMRRT} (iii) is first proved by Fomin and
Zelevinsky in ~\cite{FZ2} 3.5, 4.4 for all Dynkin diagrams
$\Gamma$ of finite type. Later we will give an analogue of the
proof in~\cite{BMRRT} for non-simply-laced finite types.

\section{Non simply-laced clusters of finite type}

\subsection{Folding cluster categories}

We call an pair $(Q,\sigma)$ {\em an admissible quiver} if $Q$ is
a (finite) quiver, and $\sigma$ an admissible automorphism of $Q$
(see~\cite{L} 12.1.1 or~\cite{DD} Example 3.5). To an admissible
quiver $(Q,\sigma)$, we associate an $\mathbb{F}_{q}$-species
$Q^{\sigma}$ (see~\cite{DD} Section 6 for detailed construction).
Conversely, each $\mathbb{F}_{q}$-species is associated to some
admissible quiver (~\cite{DD} Theorem 6.5). $\sigma$ induces a
Frobenius morphism $F$ on the algebra $\overline{\mathbb{F}}_{q}Q$
with $(\overline{\mathbb{F}}_{q}Q)^{F}=\mathbb{F}_{q}Q^{\sigma}$.
Let
$\mathcal{D}^{b}(Q)=\mathcal{D}^{b}(\overline{\mathbb{F}}_{q}Q)$
and
$\mathcal{D}^{b}(Q^{\sigma})=\mathcal{D}^{b}(\mathbb{F}_{q}Q^{\sigma})$,
then by Section~\ref{S:frobmor} we have
$\mathcal{D}^{b}(Q)^{F}\simeq \mathcal{D}^{b}(Q^{\sigma})$. Form
cluster categories
$\mathcal{C}_{Q}=\mathcal{C}_{\overline{\mathbb{F}}_{q}Q}$, and
$\mathcal{C}_{Q^{\sigma}}=\mathcal{C}_{{\mathbb{F}}_{q}Q^{\sigma}}$.
The selfequivalence $()^{[r]}:\mathcal{D}^{b}(Q)\rightarrow
\mathcal{D}^{b}(Q)$ is well-defined and commutes with the
Auslander-Reiten translation $\tau$ and the shift functor $S$, and
hence commutes with $\tau^{-1}S$. Therefore we obtain a
selfequivalence of the cluster category $\mathcal{C}_{Q}$, also
denoted by $()^{[r]}$.

\begin{Lem}\label{L:isomcluster} Let $r\in\mathbb{Z}$ and $M$ be an object in
$\mathcal{C}_{Q}$. Assume $\phi:M^{[r]}\rightarrow M$ is a
morphism in $\mathcal{C}_{Q}$, then $\phi$ is an isomorphism in
$\mathcal{C}_{Q}$ if and only if it is an isomorphism in
$\mathcal{D}^{b}(Q)$.
\end{Lem}
\begin{pf} Assume $\phi$ is an isomorphism in $\mathcal{C}_{Q}$, then $\phi$ is an isomorphism in
$\mathcal{D}^{b}(Q)$ from $M^{[r]}$ to $(\tau^{-1}S)^{i}M$ for
some $i\in\mathbb{Z}$. By a direct comparison on the dimensions of
the homology group of the two complexes we deduce that $i=0$. The
converse is obvious. $\square$
\end{pf}

A pair $(M,\phi_{M})$ is called an {\em $F$-stable object} in
$\mathcal{C}_{Q}$ if $\phi_{M}:{M}^{[1]}\rightarrow M$ is an
isomorphism in $\mathcal{C}_{Q}$, that is, an object in
$\mathcal{D}^{b}(Q)^{F}$ by Lemma~\ref{L:isomcluster}. We also say
that $M$ is an $F$-stable object. Denote by $\mathcal{C}^{F}_{Q}$
the category of $F$-stable objects whose morphisms
$f:(M,\phi_{M})\rightarrow (N,\phi_{N})$ are morphisms
$f:M\rightarrow N$ in $\mathcal{C}_{Q}$ satisfying
$f\circ\phi_{M}=\phi_{N}\circ {f}^{[1]}$.

\begin{Thm}\label{T:clustercatequiv} {\rm (i)}
$\mathcal{C}^{F}_{Q}$ is equivalent to $\mathcal{C}_{Q^{\sigma}}$
as triangulated categories. For an object $(M,\phi_{M})$ in
$\mathcal{C}^{F}(Q)$ we will denote by $M^{F}$ the corresponding
object in $\mathcal{C}_{Q^{\sigma}}$.

{\rm (ii)} Let $(M,\phi_{M})$ and $(N,\phi_{N})$ be two objects in
$\mathcal{C}^{F}_{Q}$. Then $M^{F}$ and $N^{F}$ are isomorphic in
$\mathcal{C}_{Q^{\sigma}}$ if and only if $M$ and $N$ are
isomorphic in $\mathcal{C}_{Q}$. Moreover,
$Hom_{\mathcal{C}_{Q^{\sigma}}}(M^{F},N^{F})=(Hom_{\mathcal{C}_{Q}}(M,N))^{F}$.

\end{Thm}

\begin{pf} (i)
Since $\tau^{-1}S$ commutes with $()^{[1]}$, it induces a
selfequivalence of $\mathcal{D}^{b}(Q)^{F}$ as a triangulated
category. Moreover, let $(M,\phi_{M}),(N,\phi_{N})$ be two objects
in $\mathcal{C}^{F}_{Q}$ (i.e. in $\mathcal{D}^{b}(Q)^{F}$) and
$\xi:(M,\phi_{M})\rightarrow (N,\phi_{N})$ a morphism in
$\mathcal{C}^{F}_{Q}$, i.e. $\xi:M\rightarrow N$ is a morphism in
$\mathcal{C}_{Q}$ or equivalently $\xi:M\rightarrow
\oplus_{i}(\tau^{-1}S)^{i}N$ is a morphism in $\mathcal{D}^{b}(Q)$
and the following diagram is commutative in $\mathcal{D}^{b}(Q)$.
\[\begin{CD} M^{[1]}@>\phi_{M}>> M\\
@VV\xi^{[1]}V @VV\xi V\\
\oplus(\tau^{-1}S)^{i}N^{[1]} @>\oplus(\tau^{-1}S)^{i}\phi_{N}
>> \oplus(\tau^{-1}S)^{i}N
\end{CD}\]
Therefore $\xi$ is a morphism in $\mathcal{C}_{Q}^{F}$ if and only
if it is a morphism  in $\mathcal{D}^{b}(Q)^{F}$ from
$(M,\phi_{M})$ to
$(\oplus(\tau^{-1}S)^{i}N,\oplus(\tau^{-1}S)^{i}\phi_{N})$, i.e.
it is a morphism from $(M,\phi_{M})$ to $(N,\phi_{N})$ in
$\mathcal{D}^{b}(Q)^{F}/(\tau^{-1}S)$. Thus we have proved that
$\mathcal{C}_{Q}^{F}$ is canonically equivalent to
$\mathcal{D}^{b}(Q)^{F}/\tau^{-1}S$.

Now by Theorem~\ref{T:dd} it follows that
$\mathcal{D}^{b}(Q)^{F}/(\tau^{-1}S)$ is equivalent to
$\mathcal{D}^{b}(Q^{\sigma})/(\tau^{-1}S)=\mathcal{C}^{F}_{Q}$,
and we are done. $\square$
\end{pf}

By the above proof the equivalence in
Theorem~\ref{T:clustercatequiv} and its inverse can be constructed
explicitly, but here we do not give the details. By
Lemma~\ref{L:isomcluster} the following is a consequence of the
corresponding result for the derived category
$\mathcal{D}^{b}(Q)$.

\begin{Prop}\label{P:periodicity} For an object $M$ in $\mathcal{C}_{Q}$ there exists $r\in\mathbb{N}$ such that $M^{[r]}\cong M$.
\end{Prop}

As in the derived category, for an $r$ as in
Proposition~\ref{P:periodicity}, the object $M\oplus
M^{[1]}\oplus\cdots\oplus M^{[r-1]}$ is $F$-stable. The minimal
such $r$ is called the {\em $F$-period} of $M$. We remark that the
$F$-periods of an $\overline{\mathbb{F}}_{q}Q$-module in
$mod\overline{\mathbb{F}}_{q}Q$, $\mathcal{D}^{b}(Q)$, and
$\mathcal{C}_{Q}$ all coincide.

\begin{Prop}\label{P:actionexc} Assume $\sigma$ of
order $t$. Let $M$ be an exceptional indecomposable
$\overline{\mathbb{F}}_{q}Q$-module, then the $F$-period of $M$ is
a divisor of $t$, i.e. $M^{[t]}\cong M$. Consequently, $()^{[1]}$
induces a $\sigma$-action on $ex.ind(\mathcal{C}_{Q})$, the set of
isoclasses of exceptional indecomposable objects in
$\mathcal{C}_{Q}$.
\end{Prop}
\begin{pf} Let $\{p_{j}\}_{j\in J}$ be the set of paths of $Q$,
then it is a basis of $\overline{\mathbb{F}}_{q}Q$ stable under
$F^{t}$. By Ringel~\cite{R2}, there exists a basis
$\{m_{l}\}_{l\in L}$ of $M$ such that
$p_{j}m_{l}=\sum\lambda_{jll'}m_{l'}$ where $\lambda_{jll'}=1,0$
for all $j,l,l'$. Therefore $M^{[t]}\cong M$. $\square$
\end{pf}

\begin{Prop}\label{P:indec} Let $M$ be indecomposable object in $\mathcal{C}_{Q}$ with $F$-period $r$, then $\tilde{M}=M\oplus
M^{[1]}\cdots\oplus M^{[r-1]}$ is indecomposable in
$\mathcal{C}_{Q}^{F}$. Moreover,
$D_{\tilde{M}^{F}}=\mathbb{F}_{q^{r}}$.
\end{Prop}
\begin{pf} Note that there exists $m\in\mathbb{Z}$ such that
$S^{m}M$ is an object in $mod\overline{\mathbb{F}}_{q}Q$. Since
$S^{m}$ is a selfequivalence of $\mathcal{C}_{Q}$, we have that
$S^{m}M$ is also of $F$-period $r$, and $D_{\tilde{M}^{F}}\cong
D_{(\widetilde{S^{m}M})^{F}}$. Therefore we may assume $M$ is an
object in $mod\overline{\mathbb{F}}_{q}Q$. Then
\[D_{\tilde{M}^{F}}=End_{\mathcal{C}_{Q^{\sigma}}}(\tilde{M}^{F})/radEnd_{\mathcal{C}_{Q^{\sigma}}}(\tilde{M}^{F})=End_{{\mathbb{F}}_{q}Q^{\sigma}}(\tilde{M}^{F})/radEnd_{{\mathbb{F}}_{q}Q^{\sigma}}(\tilde{M}^{F})=\mathbb{F}_{q^{r}}.\]
The last equality follows from~\cite{DD} Theorem 5.1. $\square$
\end{pf}

\begin{Thm}\label{T:foldingtilting}{\rm (i)} Let $M$ be an $F$-stable object in $\mathcal{C}_{Q}$. Then $M^{F}$ is a tilting object in $\mathcal{C}_{Q^{\sigma}}$ if and only if $M$ is  a
tilting object in $\mathcal{C}_{Q}$.

{\rm (ii)} Let $M$ be an tilting object in $\mathcal{C}_{Q}$ which
is $F$-stable, then the associated Frobenius map on $M$ induces a
Frobenius morphism on the algebra $End_{\mathcal{C}_{Q}}(M)$.
Moreover,
$End_{\mathcal{C}_{Q^{\sigma}}}(M^{F})=End_{\mathcal{C}_{Q}}(M)^{F}$
(therefore the theory introduced in Section~\ref{S:frobmor}
applies to cluster-tilted algebras). In particular $F$ induces an
admissible automorphism $\sigma_{M}$ of $Q_{M}$ such that
$Q_{M^{F}}$ is the underlying valued quiver of the
$\mathbb{F}_{q}$-species $Q_{M}^{\sigma_{M}}$.
\end{Thm}
\begin{pf} This is because
$Ext^{i}_{\mathcal{C}_{Q^{\sigma}}}(M^{F},M^{F})=Ext^{i}_{\mathcal{C}_{Q}}(M,M)^{F}$,
where $i=0,1$. $\square$
\end{pf}

\vskip10pt

Let $I=\{1,\cdots,n\}$ be the vertex set of $Q$, then
$I^{\sigma}=\{\tilde{i}|i\in I\}$ is the vertex set of
$Q^{\sigma}$ where $\tilde{i}$ is the $\sigma$-orbit of $i\in I$.
Let $T=T_{1}\oplus\cdots\oplus T_{n}$ be a tilting object in
$\mathcal{C}_{Q}$, and is $F$-stable. Let $\oplus_{k'\in
\tilde{k}}T_{k'}$ be an $F$-stable direct summand, indecomposable
in $\mathcal{C}_{Q}^{F}$. Then there is a unique set
$\{T_{k'}\}_{k'\in\tilde{k}}$ such that
$T'=\oplus_{i\not\in\tilde{k}}T_{i}\oplus
\oplus_{k'\in\tilde{k}}T_{k'}'$ is a tilting object in
$\mathcal{C}_{Q}$, and is $F$-stable. Write
$B=B_{{T^{F}}}=(b_{\tilde{i}\tilde{j}})_{\tilde{i},\tilde{j}\in
I^{\sigma}}$,
$B'=B_{{T'^{F}}}=(b_{\tilde{i}\tilde{j}}')_{\tilde{i},\tilde{j}\in
I^{\sigma}}$, $A=B_{{T}}=(a_{ij})_{i,j\in I}$, and
$A'=B_{{T'}}=(a_{ij}')_{i,j\in I}$, then
$A'=\prod_{k'\in\tilde{k}}\delta_{k'}A$, and
$B'=\delta_{\tilde{k}}B$. By Theorem~\ref{T:foldingtilting}(ii) we
have  $b_{\tilde{i}\tilde{j}}=\sum_{j'\in\tilde{j}}a_{ij'}$, and
$b_{\tilde{i}\tilde{j}}'=\sum_{j'\in\tilde{j}}a_{ij'}'$. We want a
description of $B'$ via $B$. The following proposition is proved
in~\cite{D}.
\begin{Prop}\label{P:mutationmat} $b_{\tilde{i}\tilde{j}}'=\begin{cases} -b_{\tilde{i}\tilde{j}} & {\rm \text\ if\ } \tilde{i}=\tilde{k}, {\rm \text\ or\ } \tilde{j}=\tilde{k},\\
                                                  b_{\tilde{i}\tilde{j}}+\frac{b_{\tilde{i}\tilde{k}}|b_{\tilde{k}\tilde{j}}|+|b_{\tilde{i}\tilde{k}}|b_{\tilde{k}\tilde{j}}}{2} & {\rm \text\ otherwise.}\end{cases}$

\end{Prop}
Here we give a new proof. Write
$\delta_{k}(A)=A^{(k)}=(a_{ij}^{(k)})$, then we have
\begin{Lem}\label{L:matmutationbmr}(~\cite{BMR}) $a_{ij}^{(k)}=\begin{cases} -a_{ij} & {\rm \text \ if\ } i=k {\rm \text\ or\ } j=k,\\
                             a_{ij}+\frac{a_{ik}|a_{kj}|+|a_{ik}|a_{kj}}{2} & {\rm \text\
                             otherwise.}\end{cases}$
                             \end{Lem}

From this formula we have the following
\begin{Lem}\label{L:nonadjmatmutation} Let $l_{1},\cdots,l_{t}\in\{1,\cdots,n\}$ be $t$
pairwise non adjacent vertices, i.e. $a_{l_{p}l_{p'}}=0$ for any
$p,p'\in\{1,\cdots,t\}$. Write
$\delta_{l_{t}}\cdots\delta_{l_{1}}(A)=A^{(l_{1}\cdots
l_{t})}=(a_{ij}^{(l_{1}\cdots l_{t})})$, then
\[a_{ij}^{(l_{1}\cdots l_{t})}=\begin{cases} -a_{ij} & {\rm \text\ if\ one\ of\ } i,j\in\{l_{1},\cdots,l_{t}\}, {\rm \text\ the\ other\ not,}\\
                                             a_{ij}+\sum_{p=1}^{t}\frac{a_{il_{p}}|a_{l_{p}j}|+|a_{il_{p}}|a_{l_{p}j}}{2} & {\rm\text\ otherwise.} \end{cases}\]

\end{Lem}
\begin{pf} We prove by induction on $t$.

If $t=1$, we have by Lemma~\ref{L:matmutationbmr}
\[\begin{array}{c}a_{ij}^{(l_{1})}=\begin{cases} -a_{ij} & {\rm \text \ if\ one\ of\ } i,j\in \{l_{1}\}, {\rm \text \ the\ other\ not,}\\
                             a_{ij}+\frac{a_{il_{1}}|a_{l_{1}j}|+|a_{il_{1}}|a_{l_{1}j}}{2} & {\rm \text\ otherwise.}\end{cases}
                             \end{array}\]

Suppose the statement is true for $t-1$. Then by
Lemma~\ref{L:matmutationbmr}
\[a_{ij}^{(l_{1}\cdots l_{t})}=\begin{cases} -a_{ij}^{(l_{1}\cdots l_{t-1})} & {\rm \text \ if\ one\ of\ } i,j\in \{l_{t}\}, {\rm \text \ the\ other\ not,}\\
                             a_{ij}^{(l_{1}\cdots l_{t-1})}+\frac{a_{il_{t}}^{(l_{1}\cdots l_{t-1})}|a_{l_{t}j}^{(l_{1}\cdots l_{t-1})}|+|a_{il_{t}}^{(l_{1}\cdots l_{t-1})}|a_{l_{t}j}^{(l_{1}\cdots l_{t-1})}}{2} & {\rm \text\ otherwise.}\end{cases}\]

We will use the following two conditions : (i) one of $i$,
$j\in\{l_{1},\cdots,l_{t-1}\}$, the other not; (ii) one of $i$,
$j\in\{l_{t}\}$, the other not.

Case 1 : condition (i) (ii) both hold. Namely, one of
$i,j\in\{l_{1},\cdots,l_{t-1}\}$, the other
$\not\in\{l_{1},\cdots,l_{t-1},l_{t}\}$. Then
\[\begin{array}{c}a_{ij}^{(l_{1}\cdots l_{t})}=-a_{ij}^{(l_{1}\cdots l_{t-1})}=a_{ij}=a_{ij}+\sum_{i=1}^{p}\frac{a_{il_{p}}|a_{l_{p}j}|+|a_{il_{p}}|a_{l_{p}j}}{2}.\end{array}\]

Case 2 : condition (i) holds, but condition (ii) fails. Then if
$i\in\{l_{1},\cdots,l_{t-1}\}$,
$j\not\in\{l_{1},\cdots,l_{t-1}\}$, then $a_{il_{t}}^{(l_{1}\cdots
l_{t-1})}=-a_{il_{t}}=0$; if $i\not\in\{l_{1},\cdots,l_{t-1}\}$,
$j\in\{l_{1},\cdots,l_{t-1}\}$, then $a_{l_{t}j}^{(l_{1}\cdots
l_{t-1})}=-a_{l_{t}j}=0$. Therefore
\[\begin{array}{c}a_{ij}^{(l_{1}\cdots l_{t})}=a_{ij}^{(l_{1}\cdots l_{t-1})}+\frac{a_{il_{t}}^{(l_{1}\cdots l_{t-1})}|a_{l_{t}j}^{(l_{1}\cdots l_{t-1})}|+|a_{il_{t}}^{(l_{1}\cdots l_{t-1})}|a_{l_{t}j}^{(l_{1}\cdots l_{t-1})}}{2}=-a_{ij}.\end{array}\]

Case 3 : condition (i) fails but condition (i) holds. Namely, one
of $i,j\in\{l_{t}\}$, the other
$\not\in\{l_{1},\cdots,l_{t-1},l_{t}\}$. Then
\[\begin{array}{c}a_{ij}^{(l_{1}\cdots l_{t})}=-a_{ij}^{(l_{1}\cdots l_{t-1})}=-(a_{ij}+\sum_{p=1}^{t-1}\frac{a_{il_{p}}|a_{l_{p}j}|+|a_{il_{p}}|a_{l_{p}j}}{2})=-a_{ij}.\end{array}\]
where the last equality is because either $a_{il_{p}}=0$ or
$a_{l_{p}j}=0$.

Case 4 : condition (i) (ii) both fail. Then
\[\begin{array}{c}a_{ij}^{(l_{1}\cdots l_{t})}=a_{ij}^{(l_{1}\cdots
l_{t-1})}+\frac{a_{il_{t}}^{(l_{1}\cdots
l_{t-1})}|a_{l_{t}j}^{(l_{1}\cdots
l_{t-1})}|+|a_{il_{t}}^{(l_{1}\cdots
l_{t-1})}|a_{l_{t}j}^{(l_{1}\cdots l_{t-1})}}{2}.\end{array}\] By
induction hypothesis we have
\[\begin{array}{c} a_{ij}^{(l_{1}\cdots l_{t-1})}=a_{ij}+\sum_{p=1}^{t-1}\frac{a_{il_{p}}|a_{l_{p}j|+|a_{il_{p}}|a_{l_{p}j}}}{2},\\
a_{il_{t}}^{(l_{1}\cdots
l_{t-1})}=a_{il_{t}}+\sum_{p=1}^{t-1}\frac{a_{il_{p}}|a_{l_{p}l_{t}}|+|a_{il_{p}}|a_{l_{p}l_{t}}}{2}=a_{il_{t}},\\
a_{l_{t}j}^{(l_{1}\cdots
l_{t-1})}=a_{l_{t}j}+\sum_{p=1}^{t-1}\frac{a_{l_{t}l_{p}}|a_{l_{p}j}|+|a_{l_{t}l_{p}}|a_{l_{p}j}}{2}=a_{l_{t}j}.
\end{array}\]
Therefore
\[\begin{array}{c}a_{ij}^{(l_{1}\cdots l_{t})}
=a_{ij}+\sum_{p=1}^{t-1}\frac{a_{il_{p}}|a_{l_{p}j|+|a_{il_{p}}|a_{l_{p}j}}}{2}+\frac{a_{il_{t}}|a_{l_{t}j}|+|a_{il_{t}}|a_{l_{t}j}}{2}
=a_{ij}+\sum_{p=1}^{t}\frac{a_{il_{p}}|a_{l_{p}j|+|a_{il_{p}}|a_{l_{p}j}}}{2}.
\end{array}\]

To conclude we get the desired result. $\square$

\end{pf}
$\bf{}$\\
{\bf Proof of Proposition~\ref{P:mutationmat} : } Since
$\tilde{k}$ is a $\sigma$-orbit of vertices it follows that
vertices in $\tilde{k}$ are pairwise non adjacent in $A$. The fact
that $Q_{T^{F}}$ has no loops or oriented cycles of length $2$
implies that for $m,l\in I$ all numbers in
$\{a_{m'l'}\}_{m'\in\tilde{m},l'\in\tilde{l}}$ has the same sign.
Applying Lemma~\ref{L:nonadjmatmutation} we prove the desired
result case by case.

Case 1 : $\tilde{i}=\tilde{k}$, $\tilde{j}=\tilde{k}$. Then
$b_{\tilde{i}\tilde{j}}'=\sum_{j'\in\tilde{j}}a_{ij'}'=0$.

Case 2 : $\tilde{i}=\tilde{k},\tilde{j}\neq\tilde{k}$, or
$\tilde{i}\neq\tilde{k},\tilde{j}=\tilde{k}$. Namely,
$i\in\{k_{1},\cdots,k_{s}\}$, $j\not\in\{k_{1},\cdots,k_{s}\}$ or
$i\not\in\{k_{1},\cdots,k_{s}\}$, $j\in\{k_{1},\cdots,k_{s}\}$.
Then
$b_{\tilde{i}\tilde{j}}'=\sum_{j'\in\tilde{j}}a_{ij'}'=\sum_{j'\in\tilde{j}}a_{ij'}^{(k_{1}\cdots
k_{s})}=\sum_{j'\in\tilde{j}}(-a_{ij'})=-b_{\tilde{i}\tilde{j}}$.

Case 3 : otherwise. In this case,
\[\begin{array}{c}b_{\tilde{i}\tilde{j}}'=\sum_{j'\in\tilde{j}}a_{ij'}'
=\sum_{j'\in\tilde{j}}(a_{ij'}+\sum_{k'\in\tilde{k}}\frac{a_{ik'}|a_{k'j}|+|a_{ik'}|a_{k'j}}{2})=b_{\tilde{i}\tilde{j}}+\sum_{k'\in\tilde{k}}\frac{a_{ik'}\sum_{j'\in\tilde{j}}|a_{k'j'}|+|a_{ik'}|\sum_{j'\in\tilde{j}}a_{k'j'}}{2}.\end{array}\]

For a fixed $k'\in\tilde{k}$ all numbers $a_{k'j'}$,
$j'\in\tilde{j}$ have the same sign. Therefore
\[\begin{array}{c}b_{\tilde{i}\tilde{j}}'
=b_{\tilde{i}\tilde{j}}+\sum_{k'\in\tilde{k}}\frac{a_{ik'}|\sum_{j'\in\tilde{j}}a_{k'j'}|+|a_{ik'}|\sum_{j'\in\tilde{j}}a_{k'j'}}{2}\end{array}\]
\[\begin{array}{c}=b_{\tilde{i}\tilde{j}}+\sum_{k'\in\tilde{k}}\frac{a_{ik'}|b_{\tilde{k}\tilde{j}}|+|a_{ik'}|b_{\tilde{k}\tilde{j}}}{2}
=b_{\tilde{i}\tilde{j}}+\frac{(\sum_{k'\in\tilde{k}}a_{ik'})|b_{\tilde{k}\tilde{j}}|+\sum_{k'\in\tilde{k}}|a_{ik'}|b_{\tilde{k}\tilde{j}}}{2}.\end{array}\]

We see that all numbers $a_{ik'}$, $k'\in\tilde{k}$ have the same
sign. Therefore
\[\begin{array}{c}b_{\tilde{i}\tilde{j}}'=b_{\tilde{i}\tilde{j}}+\frac{(\sum_{k'\in\tilde{k}}a_{ik'})|b_{\tilde{k}\tilde{j}}|+|\sum_{k'\in\tilde{k}}a_{ik'}|b_{\tilde{k}\tilde{j}}}{2}
=b_{\tilde{i}\tilde{j}}+\frac{b_{\tilde{i}\tilde{k}}|b_{\tilde{k}\tilde{j}}|+|b_{\tilde{i}\tilde{k}}|b_{\tilde{k}\tilde{j}}}{2}.\end{array}\]
This completes the proof. $\square$

\subsection{Folding almost positive roots}
 Folding almost positive roots is a consequence of folding
root systems, which is a strong tool in studying non-simply-laced
Kac theorem, see~\cite{DD}~\cite{DX}~\cite{Hua}~\cite{Hub}.

Let $Q$ be a Dynkin quiver and $\sigma$ an admissible
automorphism. Then the associated $\mathbb{F}_{q}$-species
$Q^{\sigma}$ is the natural $\mathbb{F}_{q}$-species of certain
valued quiver, which we also denote by $Q^{\sigma}$. Let $\Gamma$
be the underlying diagram of $Q$, then $\sigma$ induces an
automorphism of $\Gamma$, which we also denote by $\sigma$. The
associated valued graph $\Gamma^{\sigma}$ is exactly the
underlying graph of $Q^{\sigma}$. We recall that the vertex set of
$Q^{\sigma}$ (and $\Gamma^{\sigma}$) is
$I^{\sigma}=\{\tilde{i}|i\in I\}$, where $I$ is the vertex set of
$Q$ (and $\Gamma$) and $\tilde{i}$ is the $\sigma$-orbit of $i\in
I$.

$\sigma$ induces an automorphism of $\Phi(\Gamma)$ and $\Phi_{\geq
-1}(\Gamma)$ is a $\sigma$-stable subset. For
$\alpha\in\Phi(\Gamma)$ of $\sigma$-period $d$ we define
$\tilde{\alpha}=\sum_{s=0}^{d-1}\sigma^{s}\alpha$. Then
$\Phi(\Gamma)^{\sigma}=\{\tilde{\alpha}\ |\
\alpha\in\Phi(\Gamma)\}=\Phi(\Gamma^{\sigma})$ and $\Phi_{\geq
-1}(\Gamma)^{\sigma}=\{\tilde{\alpha}\ |\ \alpha\in\Phi_{\geq
-1}(\Gamma)\}=\Phi_{\geq -1}(\Gamma^{\sigma})$.

Let $I^{\sigma}=(I^{\sigma})^{+}\cup (I^{\sigma})^{-}$ be a
partition of $I^{\sigma}$ into completely disconnected subsets
then there is a partition $I=I^{+}\cup I^{-}$ of $I$ with $I^{+}$
and $I^{-}$ both completely disconnected and $\sigma$-stable, and
$(I^{\pm})^{\sigma}=(I^{\sigma})^{\pm}$. Define the compatibility
degrees $(\ ||\ )_{\Gamma}$ and $(\ ||\ )_{\Gamma^{\sigma}}$ via
these two partitions respectively, then

\begin{Lem}
$(\tilde{\alpha}||\tilde{\beta})_{\Gamma^{\sigma}}=\sum_{t=0}^{d(\beta)-1}(\alpha||\sigma^{t}\beta)_{\Gamma}$,
where $d(\beta)$ is the $\sigma$-period of $\beta$.
\end{Lem}
\begin{pf} This is because $(-\tilde{\alpha}_{i}||\tilde{\beta})_{\Gamma^{\sigma}}=\sum_{t=0}^{d(\beta)-1}(-\alpha_{i}||\sigma^{t}\beta)_{\Gamma}$ for any $i\in I$.
$\square$

\end{pf}

\subsection{Non-simply-laced cluster of finite type}

Let $Q$ be a Dynkin quiver and $\sigma$ an admissible
automorphism. By Proposition~\ref{P:actionexc} the selfequivalence
$()^{[1]}$ induces a $\sigma$-action on $ind\mathcal{C}_{Q}$. The
following lemma shows that this action commutes with taking
dimension vectors.

\begin{Lem}\label{L:comp} For any object $M$ in $\mathcal{C}_{Q}$ we have
$\underline{dim}(M^{[1]})=\sigma(\underline{dim}M)$.
\end{Lem}
\begin{pf} It suffices to prove for $M=SP_{i}$ where $P_{i}$ is the indecomposable projective corresponding to the vertex $i\in I$, and for $M$ an $\overline{\mathbb{F}}_{q}Q$-module. By the comment after~\cite{DD} Proposition4.2, we
deduce that $P_{i}^{[1]}\cong P_{\sigma(i)}$. Therefore
\[\underline{dim}((SP_{i})^{[1]})=\underline{dim}(S(P_{i})^{[1]})=\underline{dim}(SP_{\sigma(i)})=-\alpha_{\sigma(i)}=\sigma(-\alpha_{i})=\sigma(\underline{dim}(SP_{i})).\]

For $M$ an $\overline{\mathbb{F}}_{q}Q$-module, we have
\[\begin{array}{c}\underline{dim}(M^{[1]})=\sum_{i\in I}dim_{\overline{\mathbb{F}}_{q}}Hom_{Q}(P_{i},M^{[1]})\alpha_{i}
=\sum_{i\in
I}dim_{\overline{\mathbb{F}}_{q}}Hom_{Q}(P_{i}^{[-1]},M)\alpha_{i}\\
\end{array}\]\[\begin{array}{c}=\sum_{i\in
I}dim_{\overline{\mathbb{F}}_{q}}Hom_{Q}(P_{\sigma^{-1}(i)},M)\alpha_{i}
=\sum_{i\in
I}dim_{\overline{\mathbb{F}}_{q}}Hom_{Q}(P_{i},M)\alpha_{\sigma(i)}=\sigma(\underline{dim}M).\
\square
\end{array}\]
\end{pf}

\begin{Prop}\label{P:nsBMRRT1} Taking dimension vector defines a bijective correspondence between
the set $\{SP_{\tilde{i}}\ |\ \tilde{i}\in I^{\sigma}\}\cup ind\
mod(\mathbb{F}_{q}Q^{\sigma})$ and $\Phi_{\geq
-1}(\Gamma^{\sigma})$, namely, between
$ind\mathcal{C}_{Q^{\sigma}}$ and the set of cluster variables,
with $\tilde{M}_{\alpha}^{F}$ corresponding to $\tilde{\alpha}$
for any $\alpha\in\Phi_{\geq -1}(\Gamma)$.
\end{Prop}
\begin{pf} It follows from Lemma~\ref{L:comp} and Proposition~\ref{P:actionexc}. $\square$
\end{pf}

Let $\alpha,\beta\in\Phi_{\geq -1}(\Gamma)$, and ${M}_{\alpha},
M_{\beta}$ the corresponding indecomposable object of
$\mathcal{C}_{Q}$ as in Theorem~\ref{T:BMRRT}. Let $d(\alpha)$ and
$d(\beta)$ be the $\sigma$-period of $\alpha$ and $\beta$
respectively (also $F$-period of ${M}_{\alpha}$ and $M_{\beta}$
respectively). Then
\[\begin{array}{c}(\tilde{\alpha}||\tilde{\beta})_{\Gamma^{\sigma}}=\sum_{t=0}^{d(\beta)-1}(\alpha||\sigma^{t}\beta)_{\Gamma}=\sum_{t=0}^{d(\beta)-1}dim_{\overline{\mathbb{F}}_{q}}Ext^{1}_{\mathcal{C}_{Q}}({M}_{\alpha},M_{\beta}^{[1]})
=dim_{\overline{\mathbb{F}}_{q}}Ext^{1}_{\mathcal{C}_{Q}}({M}_{\alpha},\tilde{M}_{\beta}).\end{array}\]
Since
$dim_{\overline{\mathbb{F}}_{q}}Ext^{1}_{\mathcal{C}_{Q}}({M}_{\alpha},\tilde{M}_{\beta})=dim_{\overline{\mathbb{F}}_{q}}Ext^{1}_{\mathcal{C}_{Q}}({M}_{\alpha}^{[r]},\tilde{M}_{\beta}^{[r]})=dim_{\overline{\mathbb{F}}_{q}}Ext^{1}_{\mathcal{C}_{Q}}({M}_{\alpha}^{[r]},\tilde{M}_{\beta})$
for any $r\in\mathbb{Z}$, it follows that
\[\begin{array}{c}(\tilde{\alpha}||\tilde{\beta})_{\Gamma^{\sigma}}=\frac{1}{d(\alpha)}dim_{\overline{\mathbb{F}}_{q}}Ext^{1}_{\mathcal{C}_{Q}}(\tilde{M}_{\alpha},\tilde{M}_{\beta})
=\frac{1}{d(\alpha)}dim_{{\mathbb{F}}_{q}}Ext^{1}_{\mathcal{C}_{Q^{\sigma}}}(\tilde{M}_{\alpha}^{F},\tilde{M}_{\beta}^{F}).\end{array}\]
By~\ref{P:indec} we have
$D_{\tilde{M}_{\alpha}^{F}}=\mathbb{F}_{q^{d(\alpha)}}$, and hence
\[\begin{array}{c}(\tilde{\alpha}||\tilde{\beta})_{\Gamma^{\sigma}}=dim_{{\mathbb{F}}_{q^{d(\alpha)}}}Ext^{1}_{\mathcal{C}_{Q^{\sigma}}}(\tilde{M}_{\alpha}^{F},\tilde{M}_{\beta}^{F})
=dim_{D_{\tilde{M}_{\alpha}^{F}}}Ext^{1}_{\mathcal{C}_{Q^{\sigma}}}(\tilde{M}_{\alpha}^{F},\tilde{M}_{\beta}^{F}).
\end{array}\]
Therefore,

\begin{Thm}\label{T:nonsimplylaced}
{\rm (i)}
$(\tilde{\alpha}||\tilde{\beta})_{\Gamma^{\sigma}}=dim_{D_{\tilde{M}_{\alpha}^{F}}}Ext^{1}_{\mathcal{C}_{Q^{\sigma}}}(\tilde{M}_{\alpha}^{F},\tilde{M}_{\beta}^{F})$.

{\rm (ii)} The bijective correspondence in
Proposition~\ref{P:nsBMRRT1} induces a one-to-one correspondence
between isoclasses of tilting objects in
$\mathcal{C}_{Q^{\sigma}}$ and clusters.

{\rm (iii)} Two elements $\tilde{\alpha},\tilde{\beta}$ in
$\Phi_{\geq -1}(\Gamma^{\sigma})$ form an exchange pair if and
only their compatibility degree is $1$, i.e.
$(\tilde{\alpha}||\tilde{\beta})_{\Gamma^{\sigma}}=1=(\tilde{\beta}||\tilde{\alpha})_{\Gamma^{\sigma}}$.

\end{Thm}

\begin{pf} (ii) follows from (i). (iii) follows from (i) and (ii). $\square$

\end{pf}

We point out that these results are not new.
Proposition~\ref{P:nsBMRRT1} and
Theorem~\ref{T:nonsimplylaced}(i)(ii) are proved in~\cite{Z}.
Proposition~\ref{P:nsBMRRT1} is also a direct consequence of
non-simply-laced version of Gabriel's theorem (see~\cite{DR}). As
mentioned before Theorem~\ref{T:nonsimplylaced}(iii) is proved
in~\cite{FZ2}.

Let $\underline{x}=(x_{\tilde{k}})_{\tilde{k}\in I^{\sigma}}$ be a
cluster corresponding to the tilting object
$T=\oplus_{\tilde{k}\in I^{\sigma}}T_{\tilde{k}}$. Let
$\tilde{k}\in I^{\sigma}$. Let $T_{\tilde{k}}^{*}$ be the other
complement to the almost complete tilting object
$\overline{T}=\oplus_{\tilde{i}\neq \tilde{k}}T_{\tilde{i}}$. Then
$x_{\tilde{k}}'=\underline{dim}T_{\tilde{k}}^{*}$ is the unique
element in $\Phi_{\geq -1}(\Gamma^{\sigma})$ different from
$x_{\tilde{k}}$ such that
$(\underline{x}\backslash\{x_{\tilde{k}}\})\cup\{x_{\tilde{k}}'\}$
is again a cluster.

Let
$(\underline{x},B)=(\{-\tilde{\alpha}_{i}|i=1,\cdots,n\},B_{Q})$
be the initial seed and correspondingly we fix the tilting seed
$(S\mathbb{F}_{q}Q^{\sigma},B_{Q})$. For a seed
$(\underline{x}',B')=\mu_{\tilde{k}_{t}}\cdots\mu_{\tilde{k}_{1}}(\underline{x},B)$,
define
$\phi(\underline{x}',B')=\delta_{\tilde{k}_{t}}\cdots\delta_{\tilde{k}_{1}}(S\mathbb{F}_{q}Q^{\sigma},B_{Q})$.

\begin{Thm} $\phi$ is well-defined, i.e. $\phi$ does not depend on the choice of the sequence $(k_{1},\cdots,k_{t})$. Moreover, $\phi$ is a bijection from
the set of seeds to the set of tilting seeds such that for any
$\tilde{k}\in I^{\sigma}$, we have a commutative diagram
\[\UseTips\xymatrix@=15mm{(\underline{x}',B')\ar[d]^{\mu_{\tilde{k}}} \ar[r]^{\phi} & (T',B_{T'}) \ar[d]^{\delta_{\tilde{k}}}\\
(\underline{x}'',B'')\ar[r]^{\phi} & (T'',B_{T''})}\]
\end{Thm}
\begin{pf} By induction it follows from
Proposition~\ref{P:mutationmat} and
Theorem~\ref{T:nonsimplylaced}. $\square$
\end{pf}

\vskip 50pt

\bibliographystyle{amsplain}

\end{document}